\documentclass[12pt,reqno]{amsart}

\usepackage{amssymb}

\theoremstyle{plain}
\newtheorem{theorem}{Theorem}[section]
\newtheorem{proposition}[theorem]{Proposition}
\newtheorem{corollary}[theorem]{Corollary}

\theoremstyle{definition}
\newtheorem{lemma}[theorem]{Lemma}

\numberwithin{equation}{section}

\begin{document}

\title[When roots and critical points are
integers]%
{Polynomials whose roots and\\ critical points are integers}
\date{\today}
\author{Jean-Claude Evard}


\begin{abstract}

\enlargethispage{1pt}

Polynomials whose coefficients, roots, and critical \linebreak
points lie in the ring of rational integers are called \emph{nice
polynomials}. In this paper, we present a general method for
investigating such polynomials. We extend our results from the
ring of rational integers to rings of algebraic integers that are
unique factorization domains, with special interest in the ring of
Gaussian integers. We apply our method to establish strong
properties of nice polynomials whose degree is a prime power. We
present a considerable reduction of the system of equations for
nice polynomials of arbitrary degree with three roots. We
establish properties of nice antisymmetric polynomials, and
properties of the averages of the roots of the derivatives of nice
polynomials. Finally, we present new examples of nice polynomials
obtained with the help of a computer, after a considerable
simplification of the computation by our method.
\end{abstract}

\maketitle


\section{Introduction}

\noindent The problem of finding properties, characterizations,
and methods of construction of polynomials with coefficients,
roots, and critical points in the ring of rational integers is on
the list of unsolved problems published by Richard Nowakowski
\cite{Nowakowski} in \emph{The American Mathematical Monthly} in
1999. Such polynomials are called \emph{nice polynomials}. The
earliest paper on this subject was published by M. Chapple
\cite{Chapple} in 1960. The most important paper was published by
Ralph Buchholz and James MacDougall \cite{BM} in 2000. Their paper
contains a comprehensive bibliography on the subject, where we
just need to add the most recent continuation of their work
published by Eugene Flynn \cite{Flynn} in 2001, and the paper
published by Johann Walter \cite{Walter} in 1987. In the present
paper, we establish a new method to deal with nice polynomials,
and we present new properties and examples obtained with our
method. While there is still a long way to go to find a complete
description of the set of nice polynomials, it seems that many
more results can be obtained in the near future with our method.

\enlargethispage{20pt}

In Section 2, we present equivalences of nice polynomials, which
reduce the search for nice polynomials to the search of just one
representant in each equivalence class. In Section 3, we establish
the key relations to deal with nice polynomials, namely, the
system of equations~\eqref{E:3.2}, which are relations between the
roots and the critical points of a polynomial. All of the results
that follow are consequences of these relations. In Section 4, we
establish properties of the averages of the roots of the
derivatives of nice polynomials, with an application to nice
cubics. In Section 5, we deal with nice symmetric and
antisymmetric polynomials. In Section 6, we present the main
result of this paper, namely, Theorem \ref{T:6.3}, which is about
nice polynomials whose degree is a prime power, with an
application to nice quartics. In Section 7, we present an
important reduction of the system of equations \eqref{E:3.2} in
the case of nice polynomials having three roots of any
multiplicities. In Section 8, we present new examples of nice
polynomials obtained by using our method with the help of a
computer.

\noindent \textbf{Notation, definitions, and terminology.} We
denote the set of rational integers \( 0 \), \( \pm 1 \), \( \pm 2
\), \dots by $\mathbb{Z}$, the set of positive rational integers
by $\mathbb{N}$, the set of non-negative rational integers by
$\mathbb{N}_0$, the set of Gaussian integers \( a + bi \), with \(
a, b \in \mathbb{Z} \), \( i^2=-1 \), by $\mathbb{G}$, and the set
of rational numbers by $\mathbb Q$. For the sake of simplicity, we
use the term \emph{polynomial} for both, \emph{polynomial} and
\emph{polynomial function}. Let \( \frak{R} \) be a ring. We
denote by \( \frak R[x] \) the ring of polynomials in one
indeterminate $x$ with coefficients in $\frak R$. If \( \frak{R}
\) is commutative, then we say that a polynomial \( p \in
\frak{R}[x] \) of degree \( d > 0 \) \emph{splits} over $\frak R$
if $p$ has the form $p(x) = c\prod_{i=1}^d (x-x_i)$, where $c$,
$x_1$, \dots, $x_d$ belong to $\frak R$, and $c \ne 0$.

Let \( \frak{R} \) be a ring of characteristic zero possessing a
multiplicative identity element $1_{\frak R}$. Let $\Phi$ denote
the ring homomorphism from $\mathbb{Z}$ into $\frak R$ determined
by \( \Phi(1) = 1_{\frak R} \). Since \( \frak{R} \) is of
characteristic zero, the homomorphism \( \Phi \) gives a ring
isomorphism from \( \mathbb{Z} \) onto \( \Phi(\mathbb{Z}) =
\mathbb{Z} \cdot 1_{\frak R} \), so that \( \mathbb{Z} \) can be
considered as a subring of \( \frak{R} \). For the sake of
simplicity, we use the same notation $n$ for an element $n$ in
$\mathbb{Z}$ and its corresponding element $\Phi(n) = n~\cdot~
1_{\frak R}$ in $\frak R$; for example $\pm 1 = \pm 1_{\frak R}$,
$2 = 1_{\frak R} + 1_{\frak R}$, \dots.

Let $\frak D$ be an integral domain of characteristic zero. Let \(
p \in \frak{D}[x] \) be a polynomial of degree \( d > 0 \) with
coefficients \( c_0, \dots, c_d \in \frak{D} \), \( c_d \ne 0 \),
that is, \( p(x) = \sum_{k=0}^d c_kx^k \). We define the
derivative \( p' \in \frak{D}[x] \) of \( p \) as \( p'(x) =
\rule{0pt}{15pt} \sum_{k=1}^d c_k\Phi(k)x^{k-1} \), which we
simply write $p'(x) = \rule{0pt}{15pt} \sum_{k=1}^d c_kkx^{k-1}$
with the above convention. Since \( \frak{D} \) is of
characteristic zero, we have \( d \ne 0 \) in  \( \frak{D} \).
Since \( \frak{D} \) is an integral domain, it follows that \( c_d
\cdot d \ne 0 \) in  \( \frak{D} \), so that \( p' \) is of degree
\( d-1 \).

Since \( \frak{D} \) is an integral domain,  \( \frak{D} \)
possesses a (unique) field of fractions that we denote by $Q(\frak
D)$, which in turn possesses a unique algebraic closure that we
denote by $ \rule{0pt}{15pt} \overline{Q(\frak D)}$. If $d \ge 2$,
then both $p$ and $p'$ split in that algebraic closure:
\begin{equation}
p(x)= c \prod_{i=1}^d (x-x_i) \hskip20pt \mbox{and} \hskip20pt
p'(x)= c'\prod_{j=1}^{d-1} (x-x'_j), \label{E:1.1}
\end{equation}
where $c$, $c'$ are in $\frak D$, and $x_1$, \dots, $x_d$, $x'_1$,
\dots, $x'_{d-1}$ are in $\rule{0pt}{15pt} \overline{Q(\frak D)}$.
Since \( \rule{0pt}{15pt} \overline{Q(\frak D)} \) is a field, we
have that the ring $\rule{0pt}{15pt} \overline{Q(\frak D)}[x]$ of
polynomials over this field is an Euclidean domain, and therefore,
it is a unique factorization domain, which we abbreviate as UFD.
Therefore, the two decompositions \eqref{E:1.1} are unique up to a
change of the order of the factors and multiplication of the
factors by units. Consequently, the roots $x_1$, \dots, $x_d$ of
$p$ and the roots $x'_1$, \dots, $x'_{d-1}$ of $p'$ are unique in
the unique field $\rule{0pt}{15pt} \overline{Q(\frak D)}$. Because
of this, we can talk about \emph{the} roots of $p$ and \emph{the}
roots of $p'$. For the sake of simplicity, we will talk about the
$d$ roots of \( p \) and the $d-1$ roots of $p'$, even if some of
them are equal, that is, we will omit to specify \emph{repeated
according to their multiplicities}. As usual, we call the roots of
$p'$ the \emph{critical points} of $p$. We say that \( p \) is
\emph{nice} if and only if both $p$ and $p'$ split over~$\frak D$,
that is to say, iff the \( d \) roots and the \( d-1 \) critical
points of \( p \) are in $\frak D$. We say that \( p \) is
\emph{totally nice} iff \( p \), \( p' \), \dots, \( p^{(d-1)} \)
split over $\frak D$.


\section{Nicety preserving transformations}

\noindent When we establish lists of nice polynomials, we do not
want to write several times the \emph{''same''} polynomial. We say
that two nice polynomials are \emph{equivalent} if one of them can
be obtained from the other by one of the transformations described
in Proposition \ref{T:2.1} below. These transformations were
already mentioned in \cite{BM} and \cite{Caldwell} in a different
setting. In Proposition \ref{T:2.1}, we present transformations
which transform nice polynomials into nice polynomials. In
Corollary \ref{T:2.2}, we consider transformations which transform
all nice polynomials into nice polynomials, and which have an
inverse transformation with the same property.


\begin{proposition}\label{T:2.1}
{\rm [Nicety preserving transformations].} Let $\frak D$ be an
integral domain of characteristic zero. Let $p \in \frak D[x] $ be
a nice polynomial of degree \( d \ge 2 \). Then we obtain another
nice polynomial $q \in \frak D[x]$ of degree \( d \) from $p$ in
the following three cases: \\
{\rm \textbf{(a)}} $q(x) = ap(ux-b)$, when $a, b, u \in \frak D$,
$a \ne
0$, and $u$ is a unit. \\
{\rm \textbf{(b)}} $q(x) = a^dp(a^{-1}x)$, when $a$ is a nonzero
element of $\frak D$, and $a^{-1}$ denotes its multiplicative
inverse in
$Q(\frak D)$.\\
{\rm \textbf{(c)}} $q(x) = p(ax)$, when $a \in \frak D$ is a
common divisor of the roots and critical points of $p$.
\end{proposition}

\begin{proof}
Since \( p \) is nice, we have that $p$ and $p'$ split over $\frak
D$, that is
\begin{equation}
p(x) = c \prod_{i=1}^d (x-x_i) \hskip20pt \mbox{and} \hskip20pt
p'(x)= c'\prod_{j=1}^{d-1} (x-x'_j), \notag
\end{equation}
where $c$, $c'$, $x_1$, \dots, $x_d$, $x'_1$, \dots, $x'_{d-1}$
are in $\frak D$.

\noindent \textbf{(a)} Since $u$ is unit of $\frak D$, $u$
possesses a multiplicative inverse $u^{-1} \in \frak D$. Then
\begin{align}
q(x) &= ap(ux-b) = ac\prod_{i=1}^d \big[(ux-b)-x_i\big]
\notag \\
&= ac\prod_{i=1}^d \Big[ u \big[ x-u^{-1}(b + x_i) \big] \Big] =
acu^d\prod_{i=1}^d \big[ x-u^{-1}(b + x_i) \big]. \notag
\end{align}
By the well-known properties of the algebraic derivative, we get
\begin{align}
q'(x) &= ap'(ux - b)u = auc'\prod_{j=1}^{d-1}
\big[(ux-b)-x'_j\big] \notag \\
&= auc'\prod_{j=1}^{d-1} \Big[ u \big[ x-u^{-1}(b + x'_j) \big]
\Big] = au^dc'\prod_{j=1}^{d-1} \big[ x-u^{-1}(b + x'_j) \big].
\notag
\end{align}
Thus both $q$ and $q'$ split over $\frak D$, that is to say, \( q
\) is nice.

\noindent \textbf{(b)} We have
\begin{align}
q(x) &= a^dc\prod_{i=1}^d (a^{-1}x-x_i) = a^dc\prod_{i=1}^d
\big[a^{-1}(x-ax_i)\big] = c\prod_{i=1}^d (x-ax_i), \notag
\end{align}
which shows that $q \in \frak D[x]$, and $q$ splits over $\frak
D$. By the well-known properties of the algebraic derivative, we
get
\begin{align}
q'(x) &= a^dp'(a^{-1}x)a^{-1} = a^{d-1}c'\prod_{j=1}^{d-1}
(a^{-1}x-x'_j) \notag \\
&= a^{d-1}c'\prod_{j=1}^{d-1} \big[ a^{-1}(x-ax'_j) \big] =
c'\prod_{j=1}^{d-1} (x-ax'_j).  \notag
\end{align}
Thus $q' \in \frak D[x]$, and $q'$ splits over $\frak D$.

\noindent \textbf{(c)} By hypothesis, there exist $y_1$, \dots,
$y_d$, $y'_1$, \dots, $y'_{d-1}$ in $\frak D$ such that $x_1 =
ay_1$, \dots, $x_d = ay_d$, $x'_1 = ay'_1$, \dots, $x'_{d-1} =
ay'_{d-1}$. Consequently,
\begin{equation}
q(x) = c\prod_{i=1}^d (ax-x_i) = c\prod_{i=1}^d \big[a(x-y_i)\big]
= a^dc\prod_{i=1}^d (x-y_i). \notag
\end{equation}
By the well-known properties of the algebraic derivative, we get
\begin{equation}
q'(x) = p'(ax)a = ac'\prod_{j=1}^{d-1} (ax-x'_j) =
ac'\prod_{j=1}^{d-1} \big[a(x-y'_j)\big] = a^dc'\prod_{j=1}^{d-1}
(x-y'_j). \notag
\end{equation}
Thus both $q$ and $q'$ split over $\frak D$.
\renewcommand{\qedsymbol}{$\blacksquare$}
\end{proof}

\begin{corollary}\label{T:2.2}
{\rm [Nicety preserving transformations whose inverses are also
nicety preserving].} Let $\frak D$ be an integral domain of
characteristic zero, let $p \in \frak D[x] $ be a polynomial of
degree \( d \ge 2 \), let $u_1$ and $u_2$ be units of $\frak D$,
let $a \in \frak D$, and let $q(x) = u_1 p(u_2 x + a)$. Then \( p
\) is nice iff \( q \) is nice. In particular, this equivalence
holds when $q(x) = p(-x)$.
\end{corollary}

\begin{proof}
By hypothesis, $u_1$ and $u_2$ possess multiplicative inverses
$u_1^{-1}$ and $u_2^{-1}$ in $\frak D$. Let $y = u_2 x + a$. Then
$x = u_2^{-1}(y - a)$, so that $q(x) = u_1 p(u_2 x + a)$ gives
$p(y) = u_1^{-1}q(u_2^{-1}y - u_2^{-1} a)$, and the conclusion
follows by Proposition \ref{T:2.1} (a).
\renewcommand{\qedsymbol}{$\blacksquare$}
\end{proof}


\section{Relations between the roots and the\\
critical points of a polynomial}

\noindent In Proposition \ref{T:3.1}, we recall the well-known
relations between the roots and the coefficients of a polynomial.
In Corollary \ref{T:3.2}, we deduce a system of relations between
the roots and the critical points of a polynomial, namely,
Equations \eqref{E:3.2}, which is the key tool for dealing with
nice polynomials.


\noindent \textbf{Notation.} Let $\frak R$ be a commutative ring.
For every positive integer~$k$, we define the $k^\text{th}$
\emph{elementary symmetric polynomial~$s_k$ in the~$m$
indeterminates} $ x_1, \dots, x_m$ by
\begin{displaymath}
s_k(x_1, \dots, x_m) = \sum_{1 \le i_1 <i_2< \cdots <i_k \le
m}x_{i_1} \cdots x_{i_k},
\end{displaymath}
and we define $s_0(x_1, \dots, x_m) =1$.

Let us recall the following well-known proposition, which is due
to Fran\c{c}ois Vi\`{e}te (1540--1603) for polynomials of low
degrees, and to Albert Girard (1595--1632) for the general case
\cite[Equation (1.2.3) p. 6, and p. 63]{RS}.

\begin{proposition}\label{T:3.1}
{\rm [Expression of the coefficients of a polynomial in terms of
its roots].} Let $\frak R$ be a commutative ring. Let $p \in \frak
R[x]$ be a polynomial of nonzero degree $n$ which splits over
$\frak R$, that is
\begin{equation}
p(x) = \sum_{k=0}^n c_kx^k  = c_n\prod_{i=1}^n (x-x_i), \notag
\end{equation}
where $c_0$, \dots, $c_n$, $x_1$, \dots, $x_n$ are in $\frak R$,
and \( c_n \ne 0 \). Then for every $k \in \{0,\dots,n\}$, we have
\begin{equation}
c_k = (-1)^{n-k}c_ns_{n-k}(x_1,\dots,x_n). \notag
\end{equation}
\end{proposition}

\begin{corollary}\label{T:3.2}
{\rm [Relations between the roots and the critical points of a
polynomial].} Let $\frak D$ be an integral domain of
characteristic zero. Let $p$ and $q$ $\in \frak D[x]$ be
polynomials of degrees $d \ge 2$ and $d-1$, respectively, that
split over $\frak D$, that is
\begin{alignat}{2}
p(x)&= \sum_{k=0}^d c_kx^k && = c_d \prod_{i=1}^d (x-x_i),  \notag\\
q(x)&= \sum_{k=0}^{d-1} c'_kx^k && = c'_{d-1} \prod_{i=1}^{d-1}
(x-x'_i), \label{E:3.1}
\end{alignat}
where $c_0$, \dots, $c_d$, $c'_0$, \dots, $c'_{d-1}$, $x_1$,
\dots, $x_d$, $x'_1$, \dots, $x'_{d-1}$ are in $\frak D$, and $c_d
\ne 0$, $c'_{d-1} \ne 0$. Then $q$ is the derivative of $p$ if and
only if $c'_{d-1}=dc_d$, and for every $k \in \{1,\dots,d-1\}$,
\begin{equation}
(d-k)s_k(x_1,\dots,x_d) = ds_k(x'_1,\dots,x'_{d-1}). \label{E:3.2}
\end{equation}
\end{corollary}

\begin{proof}
By substituting $k = j-1$ into \eqref{E:3.1}, we get that the
relation $q = p'$ is equivalent to $\sum_{j=1}^d c'_{j-1}x^{j-1} =
\sum_{j=1}^d c_jjx^{j-1}$, which implies that
\begin{equation}
c'_{j-1} = jc_j, \hskip30pt \forall j \in \{1, \dots, d \}.
\label{E:3.3}
\end{equation}
For every $k \in \{0, \dots, d-1\}$, let
\begin{equation}
S_k = s_k(x_1, \dots, x_d) \hskip30pt \mbox{ and } \hskip30pt S'_k
= s_k(x'_1, \dots, x'_{d-1}). \notag
\end{equation}
By applying Proposition \ref{T:3.1} to both sides of Equation
\eqref{E:3.3}, observing that \( (d-1) - (j-1) = (d-j) \), we
obtain that $q = p'$ if and only if
\begin{equation}
(-1)^{d-j}c'_{d-1}S'_{d-j} = j(-1)^{d-j}c_dS_{d-j}, \hskip30pt
\forall j \in \{1, \dots, d \}, \notag
\end{equation}
that is, with $k=d-j$,
\begin{equation}
c'_{d-1}S'_k = (d-k)c_dS_k, \hskip30pt \forall k \in \{0, \dots,
d-1 \}, \notag
\end{equation}
that is, when $k=0$, $c'_{d-1} = dc_d$, and otherwise
\begin{equation}
dc_dS'_k = (d-k)c_dS_k, \hskip30pt \forall k \in \{1, \dots, d-1
\}. \notag
\end{equation}
Since $c_d \ne 0$ and $\frak D$ is an integral domain, we can
simplify both sides of this relation by $c_d \ne 0$, which gives
the conclusion.
\renewcommand{\qedsymbol}{$\blacksquare$}
\end{proof}


\section{Averages of the roots}

\noindent

In Proposition \ref{T:4.1}, we present properties and relations
between the averages of the roots of polynomials and their
derivatives of all orders, and in Corollary \ref{T:4.2}, we
present an application to nice cubics.

\noindent \textbf{Definition.} Let $\frak D$ be an integral domain
of characteristic zero. Let \( p \in \frak D[x] \) be a polynomial
of degree $n \ge 1$. Then \( p \) splits over \( \rule{0pt}{15pt}
\overline{Q(\frak D)} \), that is, $p(x)= c \prod_{i=1}^n
(x-x_i)$, where $c \in \frak{D}$, and $x_1$, \dots, $x_n$ are in
\( \rule{0pt}{15pt} \overline{Q(\frak D)} \). Since \( \frak{D} \)
is of characteristic zero, we have that $n \ne 0$ in $\frak D$,
which allows us to define the \emph{average \( \mathcal A(p)\) of
the roots of} $p$, as
\begin{equation}
\mathcal A(p) = \frac{x_1 + x_2+ \cdots + x_d}{n} \in
\overline{Q(\frak D)}. \notag
\end{equation}

Another version of the following proposition, valid for totally
nice polynomials, can be found in \cite{Carroll}.

\begin{proposition} \label{T:4.1} 
{\rm [The average of the roots of any polynomial is equal to the
average of the roots of each of its derivatives, this common
average is rational even if some of the roots are irrational, and
if the polynomial is nice, then this common average is an
integer.]} Let \( \frak D \) be an integral domain of
characteristic zero. Let \( p \in \frak{D}[x] \) be a polynomial
of degree \( d \ge 2 \). Let \( r \in \rule{0pt}{15pt}
\overline{Q(\frak D)} \) denote the only root of \( p^{(d-1)} \).
Then
\begin{displaymath}
\mathcal A(p) = \mathcal A(p') = \mathcal A(p'') = \cdots =
\mathcal A(p^{(d-1)}) = r \in Q(\frak D).
\end{displaymath}
Furthermore, if \( \frak{D} \) is a UFD and \( p \) is nice, then
\( r \in \frak{D} \).
\end{proposition}

\begin{proof}
Since \( p \) and \( p' \) split in \( \rule{0pt}{15pt}
\overline{Q(\frak D)} \), they have the form
\begin{equation}
p(x) = c\prod_{i=1}^d (x-x_i), \hskip40pt p'(x) =
c'\prod_{i=1}^{d-1} (x-x'_i), \notag
\end{equation}
where \( c, c' \in \frak{D} \), and \( x_1, \dots, x_d, x'_1,
\dots, x'_{d-1} \in \overline{Q(\frak D)} \).  Let
\begin{equation}
S_1 = s_1(x_1,\dots,x_d) \hskip25pt \mbox{and} \hskip25pt S'_1 =
s_1(x'_1,\dots,x'_{d-1}). \notag
\end{equation}
By Equation \eqref{E:3.2}, we have \( (d-1)S_1 = dS'_1 \), which
gives \( \frac{S_1}{d} = \frac{S'_1}{d-1} \), that is \(
\mathcal{A}(p) = \mathcal{A}(p') \). Thus, the averages of the
roots of a polynomials and the roots of its first derivative are
equal. It follows by induction that the averages of the roots of
all of the derivatives of \( p \) are equal. Since the polynomial
\( p^{(d-1)} \in \frak{D}[x] \) is of degree one, it has the form
\( p(x) = ax + b \), with \( a, b \in \frak{D} \) and \( a \ne 0
\). Consequently, \( \mathcal{A} (p^{(d-1)}) = r = -\frac{b}{a}
\in Q(\frak{D}) \).

\enlargethispage{20pt}

Assume that \( \frak{D} \) is a UFD and \( p \) is nice. Then \(
x_1 \), \dots, \( x_d \), \( x'_1 \), \dots, \( x'_{d-1} \) are in
\( \frak{D} \), and so are \( S_1 \) and \( S'_1 \). Since every
common divisor of $d$ and $(d-1)$ in~$\frak D$ divides the
difference $d-(d-1) = 1$, we have that \( d \) and \( (d-1) \) do
not have any common prime divisors in \( \frak{D} \). This
together with the equation $(d-1)S_1 = dS'_1$ in \( \frak{D} \)
imply that all of the prime factors of the prime decomposition of
\( d \) are in the prime decomposition of~\( S_1 \). Consequently,
$d$ divides $S_1$, so that \( r = \mathcal{A}(p) = \frac{S_1}{d}
\in \frak{D} \).
\renewcommand{\qedsymbol}{$\blacksquare$}
\end{proof}

\begin{corollary}\label{T:4.2}
Let \( \frak D \) be a UFD of characteristic zero. Then every nice
cubic over \( \frak{D} \) is totally nice.
\end{corollary}

\begin{proof}
Let \( p \) be a nice cubic over \( \frak{D} \). Then by
definition, \( p \) and \( p' \) split over \( \frak{D} \).
Furthermore, by Proposition \ref{T:4.1}, the only root of \( p''
\) is in \( \frak{D} \), that is to say, \( p'' \) also splits
over \( \frak{D} \), so that \( p \) is totally nice.
\renewcommand{\qedsymbol}{$\blacksquare$}
\end{proof}


\section{Nice symmetric and antisymmetric polynomials}

\noindent This section deals with symmetric and antisymmetric nice
polynomials. It is restricted to polynomials with roots in \(
\mathbb{Z} \), because we use Rolle's theorem, considering such
polynomials as real-valued functions of a real variables. For
example, if \( p \) is a symmetric polynomial of degree \( d \ge 2
\), with real roots, and center \( C \) of the set of its roots,
and if \( C \) is not a root of \( p \), then \( p' \) is
antisymmetric with center \( C \), consequently, \( C \) is a root
of \( p' \), and the multiplicity of this root is one because of
Rolle's theorem and the fact that \( p' \) is of degree \( (d-1)
\). This is not true for a symmetric polynomial with complex roots
(See Example B4 in Section 8). Theorem \ref{T:5.1} gives necessary
conditions for the existence of nice antisymmetric polynomials.
Corollary \ref{T:5.2} deals with the special case where the center
of a nice antisymmetric polynomial is a root of multiplicity one
of this polynomial. Corollary \ref{T:5.3} shows that such
polynomials cannot be totally nice.

\enlargethispage{10pt}

\noindent \textbf{Definitions.} Let \( \frak{R} \) be a ring, and
let \( p \in \frak{R}[x] \) be a non-constant polynomial. We say
that \( p \) is \emph{symmetric} iff there exists \( C \in
\frak{R} \) such that \( p(C - a) = p(C + a) \) for all \( a \in
\frak{R} \). We call \( C \) the \emph{center} of~\( p \). We say
that~\( p \) is \emph{antisymmetric} iff there exists \( C \in
\frak{R} \) such that \( p(C - a) = -p(C + a) \) for all \( a \in
\frak{R} \). We call \( C \) the \emph{center} of \( p \).

\begin{theorem}\label{T:5.1}
Let \( p \in \mathbb{Z}[x] \) be a nice antisymmetric polynomial
of degree \( d \ge 3 \) and center \( C \) with more than one
root. Let \( m_0 \) denote the multiplicity of the root \( C \) of
\( p \). Let \( g = \gcd(m_0,d) \). Then the integers \(
\frac{m_0}{g} \) and \( \frac{d}{g} \) are odd squares.
\end{theorem}

\begin{proof}
Let \( q(x) = p(x+C) \). Then \( q \) is a nice over \( \mathbb{Z}
\) by Proposition~\ref{T:2.1}~(a), and it is easy to check that \(
q \) is antisymmetric with center zero. Consequently, \( r \in
\mathbb{Z} \) is a root of \( q \) iff \( (-r) \) is also a root
of \( q \) with same multiplicity. Since \( C \) is a root of \( p
\) of multiplicity \( m_0 \), we have that zero is a root of \( q
\) of multiplicity \( m_0 \). Since \( p \) has more than one
root, we have that the number \( N \) of positive roots of \( q \)
is greater than zero. It follows that \( q \) has the form
\begin{displaymath}
q(x) = c x^{m_0} \prod_{i=1}^N (x^2-x_i^2)^{m_i},
\end{displaymath}
where \( c \in \mathbb{Z} \), \( x_1, \dots, x_N \in \mathbb{N} \)
are distinct, and \( m_1, \dots, m_N \) are such that \( d = m_0 +
2m_1 + \cdots + 2m_N \). Since \( q \) is an odd function, we know
from Taylor polynomials that all of the terms of \( q \) are of
odd degrees, which implies that \( d \) and \( m_0 \) are odd.
Furthermore, all of the terms of \( q' \) are of even degree, so
that \( q' \) is symmetric with center zero. Consequently, \( r'
\in \mathbb{Z} \) is a root of \( q' \) iff \( (-r') \) is also a
root of \( q' \) with the same multiplicity. It is well known from
Algebra that each root of~\( q \) of multiplicity \( m \) gives a
root of \( q' \) of multiplicity at least \( (m-1) \). In
addition, since \( q \) has \( (2N + 1) \) distinct roots in \(
\mathbb{Z} \), we have by virtue of Rolle's theorem that \( q' \)
has \( 2N \) additional roots. Taking into account that the degree
of \( q' \) is \( (d-1) \), we obtain that \( q' \) has the form
\begin{displaymath}
q'(x) = c d x^{m_0-1} \prod_{i=1}^N (x^2-x_i^2)^{m_i-1}
\prod_{j=1}^N (x^2 - {x'_j}^2),
\end{displaymath}
where \( x'_1, \dots, x'_N \in \mathbb{N} \) are distinct. By
applying Equation \eqref{E:3.2} with \( k = d-m_0 \), we obtain
\begin{displaymath}
m_0 \cdot \prod_{i=1}^N x_i^{m_i} \cdot (-x_i)^{m_i} = d \cdot
\prod_{i=1}^N x_i^{m_i-1} \cdot (-x_i)^{m_i-1} \cdot \prod_{j=1}^N
x'_i \cdot (-x_i).
\end{displaymath}
Let \( \bar{m} = \frac{m_0}{g} \) and \( \bar{d} = \frac{d}{g} \).
Then the above equation becomes
\begin{displaymath}
\bar{m} \left(\prod_{i=1}^N x_i^{m_i} \right)^2 = \bar{d} \left(
\prod_{i=1}^N x_i^{m_i-1} \cdot \prod_{j=1}^N x'_i \right)^2.
\end{displaymath}
Since \( \bar{m} \) and \( \bar{d} \) have no common prime factor,
the above relation implies that all of the prime factors of the
prime decompositions of \( \bar{m} \) and \( \bar{d} \) are raised
to even powers, which implies that \( \bar{m} \) and \( \bar{d} \)
are squares. Since \( m_0 = g \bar{m} \) and \( d = g \bar{d} \)
are odd, we have that \( \bar{m} \) and \( \bar{d} \) are odd too.
\renewcommand{\qedsymbol}{$\blacksquare$}
\end{proof}

\begin{corollary}\label{T:5.2}
{\rm \textbf{(a)}} Let \( p \in \mathbb{Z}[x] \) be a nice
antisymmetric polynomial of degree \( d \ge 3 \) whose center \( C
\in \mathbb{Z} \) is a root of \( p \) of multiplicity one. Then
\( d \) is an odd square.
\\
{\rm \textbf{(b)}} Let \( p \in \mathbb{Z}[x] \) be a symmetric
polynomial of degree \( d \ge 4 \) whose center \( C \in
\mathbb{Z} \) is not a root of \( p \), and whose derivative is
nice. Then \( (d-1) \) is an odd square.
\end{corollary}

\begin{proof}
\textbf{(a)} This is a direct consequence of Theorem \ref{T:5.1}
applied with \( m_0 = 1 \).

\noindent \textbf{(b)} By using that the derivative of a symmetric
polynomial is antisymmetric, using the hypothesis that \( C \) is
not a root of \( p \), using that the center is a root of all
antisymmetric polynomials, using that the degree of \( p' \) is \(
d-1 \), and using Rolle's theorem, we deduce that \( p' \) is
antisymmetric, and its center is a root of \( p' \) of
multiplicity one, which implies (b) by~(a).
\renewcommand{\qedsymbol}{$\blacksquare$}
\end{proof}

\noindent \textbf{Existence of non-trivial nice antisymmetric
polynomials:} In Example A 2 of Section 8, we provide examples of
nice antisymmetric polynomials of degree \( d \), as described in
Corollary \ref{T:5.2} (a), for every odd square \( d \ge 9 \). In
Example A 3, we provide examples of non-trivial nice symmetric
polynomials of degree \( n^2+1 \) whose derivative is also nice
for every odd positive integer \( n \).

\begin{corollary}\label{T:5.3}
{\rm \textbf{(a)}} Let \( p \in \mathbb{Z}[x] \) be a nice
antisymmetric polynomial of degree \( d \ge 5 \) whose center \( C
\) is a root of \( p \) of multiplicity one. Then \( p'' \) is not
nice.
\\
{\rm \textbf{(b)}} Let \( p \in \mathbb{Z}[x] \) be a symmetric
polynomial of degree \( d \ge 6 \) whose center \( C \) is not a
root of \( p \), and whose derivative is nice. Then \( p''' \) is
not nice.
\end{corollary}

\begin{proof}
\noindent \textbf{(a)} Since \( p \) is antisymmetric of degree \(
d \) with center \( C \), we have that \( p' \) is symmetric of
degree \( d-1 \) with same center, and \( p'' \) is antisymmetric
of degree \( d-2 \) with same center. It is easy to check using
Rolle's theorem that the hypothesis implies that \( C \) is not a
root of \( p' \), and  \( C \) is a root of \( p'' \) of
multiplicity one. If \( p'' \) were nice, then both \( d \) and \(
d-2 \) would be squares by Corollary \ref{T:5.2} (a), which is
impossible, because the distance between two nonzero distinct
squares is at least three. Therefore, \( p'' \) is not nice.

\noindent \textbf{(b)} The hypothesis implies that \( p' \) is a
nice antisymmetric polynomial of degree \( d \ge 5 \) whose center
\( C \in \mathbb{Z} \) is a root of \( p \) of multiplicity one.
This implies the conclusion by (a).
\renewcommand{\qedsymbol}{$\blacksquare$}
\end{proof}


\section{Nice polynomials whose degree is a prime-power}

\noindent In this section, we present the main result of this
paper, which is about nice polynomials whose degree is a
prime-power. This result holds for polynomials whose roots and
critical points belong to a UFD \( \frak D \) of characteristic
zero such that \( \frak D \cap \mathbb Q = \mathbb Z \) in \(
Q(\frak D) \). Proposition~\ref{T:6.1} shows that this condition
is satisfied by every free integral domain, and in particular, by
every ring of algebraic integers of an algebraic number field.
Lemma \ref{T:6.2} is the key lemma to extend our main theorem from
the ring \( \mathbb Z \) of rational integers to rings of
algebraic integers that are UFDs. Theorem \ref{T:6.3} about nice
polynomials of degree \( d = b^e \), where \( b \) is a prime, is
the main result of this paper, with a stronger conclusion when the
base \( b \) divides the exponent \( e \). Corollary \ref{T:6.4}
presents an important application to nice quartics.

\noindent \textbf{Definitions.} Let \( \frak D \) be an integral
domain. We say that a subset~\( \mathcal B \) of~\( \frak D \) is
a \( \mathbb Z \)-\emph{basis} of \( \frak D \) iff every element
\( a \) of  \( \frak D \) can be expressed in a unique way as a
finite \( \mathbb Z \)-linear combination  of elements of \(
\mathcal B \), that is, iff for every \( a \in \frak D \), there
exists a unique family \( (a_\beta)_{\beta \in \mathcal B} \) of
elements of \( \mathbb Z \), with only a finite number of nonzero
elements, such that \( a = \sum_{\beta \in \mathcal B} a_\beta
\beta \). We say that \( \frak D \) is \emph{free} iff \( \frak D
\) possesses a \( \mathbb Z \)-\emph{basis}.

Let \( \frak{D} \) be an integral domain of characteristic zero.
Then the ring homomorphism \( \Phi \) from \( \mathbb{Z} \) into
\( \frak{D} \) determined by \( \Phi(1) = 1_{\frak{R}} \) can be
extended to a field homomorphism from \( \mathbb{Q} \) to \(
Q(\frak{D}) \). Furthermore, since  \( \frak{D} \) is of
characteristic zero, the homomorphism \( \Phi \) gives a ring
isomorphism from \( \mathbb{Z} \) onto \( \Phi(\mathbb{Z}) =
\mathbb{Z} \cdot 1_{\frak R} \), and a field isomorphism from \(
\mathbb{Q} \) onto \( \Phi(\mathbb{Q}) = \mathbb{Q} \cdot 1_{\frak
R} \). Consequently, \( \mathbb{Z} \) can be considered as a
subring of \( \frak{D} \), \( \mathbb{Q} \) can be considered as a
subfield of \( Q(\frak{D}) \), and both \( \mathbb{Z} \) and \(
\mathbb{Q} \) can be considered as subrings of \( Q(\frak{D}) \).

\begin{proposition}\label{T:6.1}
Every free integral domain \( \frak D \) is a ring of
characteristic zero satisfying the condition \( \frak D \cap
\mathbb Q = \mathbb Z \) in \( Q(\frak D) \). In particular, every
ring of algebraic integers of an algebraic number field has these
two properties.
\end{proposition}

\begin{proof}
Let \( \frak D \) be an integral domain with a \( \mathbb Z
\)-basis \( \mathcal B \). Assume that \( \frak D \) has nonzero
characteristic \( n \in \mathbb N \). Then \( n a = 0 \) for all
in \( a \in \frak D \). Since \( 1 \ne 0 \) in the integral domain
\( \frak D \), the basis \( \mathcal B \) has at least one element
\( \beta_0 \). Then \( \beta_0 = 1 \beta_0 = (1+n)\beta_0 \), so
that \( \beta_0 \) can be expressed in two different ways as a
finite \( \mathbb Z \)-linear combination of elements of the basis
\( \mathcal B \). This contradicts the definition of a basis.
Therefore \( \frak D \) has characteristic zero.

Clearly \( \mathbb Z \subseteq \frak D \cap \mathbb Q \). Let us
prove that \( \frak D \cap \mathbb Q \subseteq \mathbb Z \). Let
\( c \in \frak D \cap \mathbb Q \). Then \linebreak there exist \(
a, b \in \mathbb Z \) with \( b \ne 0 \) such that \( c =
\frac{a}{b} \). By the expression of~\( (c \beta_0) \) in the \(
\mathbb Z \)-basis \( \mathcal B \), there exists a unique family
\( (n_\beta)_{\beta \in \mathcal B} \) of elements of \( \mathbb Z
\), with only a finite number of nonzero elements, such that \( c
\beta_0 = \sum_{\beta \in \mathcal B} n_\beta \beta \). This gives
\( a \beta_0 = b c \beta_0 = \sum_{\beta \in \mathcal B} b n_\beta
\beta \). Since \( a \), \( b \), and \( n_{\beta_0} \) are in \(
\mathbb{Z} \), it follows that \( a = b n_{\beta_0} \), that is \(
bc = b n_{\beta_0} \). Since~\( \frak D \) is of characteristic
zero and \( b \ne 0 \) in \( \mathbb Z \), we have that \( b \ne 0
\) in \( \mathcal D \). Consequently, \( bc = b n_{\beta_0} \)
implies that \( c = n_{\beta_0} \in \mathbb Z \). Thus, \( \frak D
\cap \mathbb Q = \mathbb Z \).

Let \( \frak F \) be an algebraic number field. Let \( \mathcal
D_{\frak F} \) denote the ring of algebraic integers of  \( \frak
F \). Then \( \mathcal D_{\frak F} \) is an integral domain,
because it is a subring of the field \( \frak F \). Furthermore,
\( \mathcal D_{\frak F} \) is free by \cite[Theorem 1.69,
p.39]{Mollin}, and the conclusion follows by the first part of
this proof.
\renewcommand{\qedsymbol}{$\blacksquare$}
\end{proof}

\begin{lemma}\label{T:6.2}
{\rm [Key lemma for the extension of the main theorem from the
ring \( \mathbb Z \) to rings of algebraic integers].} Let $\frak
D$ be an integral domain of characteristic zero satisfying \(
\frak D \cap \mathbb Q = \mathbb Z \). Let $a, b, e \in
\mathbb{N}$ be such that \( a < b^e \). Let \( h \in \mathbb N_0
\) denote the largest exponent such that \( b^h \) divides \( a \)
in $\frak D$. Then \( h < e \).
\end{lemma}

\begin{proof}
By definition of \( h \), there exists \( c \in \frak D \) such
that \( a = b^h c \) in $\frak D$. Then \( c = \frac{a}{b^h} \in
\frak D \cap \mathbb Q = \mathbb{Z} \), so that \( c \in \mathbb Z
\). Since the ring \( \frak D \) is of characteristic zero, the
relation \( a = b^h c \) in \( \mathbb Z \cdot 1_{\frak D} \)
implies that the same relation holds in~\( \mathbb Z \). Since \(
a \) and \( b \) are positive integers, the relation \( a = b^h c
\) in \( \mathbb Z \) implies that the integer \( c \) is
positive, which implies that \( c \ge 1 \). Consequently, we have
\( b^h \le b^h c = a < b^e\) in \( \mathbb{N} \), which implies
the conclusion.
\renewcommand{\qedsymbol}{$\blacksquare$}
\end{proof}

\begin{theorem}\label{T:6.3}
{\rm [Main theorem: Properties of nice polynomials whose degree is
a power of a prime].} Let $\frak D$ be a UFD of characteristic
zero satisfying \( \frak D \cap \mathbb Q = \mathbb Z \). Let \(
b, e \in \mathbb{N} \) be such that \( b \) is prime in \( \frak D
\). Let \( d = b^e \). Let \( p \in \frak{D}[x] \) be a nice
polynomial of degree \( d \) with a root at zero. Then the roots
and critical points of \( p \) have the following properties: \\
{\rm \textbf{(a)}} All of the roots of \( p \) are multiples of \(
b \). \\
{\rm \textbf{(b)}} At most \( e \) of the critical points of \( p
\) are not multiples of \( b \). \\
{\rm \textbf{(c)}} In the special case where zero is a root of \(
p \) of multiplicity one, and~\( b \) divides \( e \) in \( \frak
D \), we have the following additional two properties: \\
{\rm \textbf{(c1)}} At least one nonzero root of \( p \) is a
multiple of \( b^2 \). \\
{\rm \textbf{(c2)}} At most \( e-1 \) of the critical points of \(
p \) are not multiples of \( b \).
\end{theorem}

\begin{proof}
By hypothesis, $p$ and $p'$ have the form
\begin{displaymath}
p(x)= cx \prod_{i=1}^{d-1} (x-x_i),  \hskip30pt p'(x)= c'
\prod_{j=1}^{d-1} (x-x'_j),
\end{displaymath}
where $c, c' \in \frak D$, the roots $x_1, \dots, x_{d-1} \in
\frak D$ of \( p \) are numbered in such a way that there exists
\( m \in \{0, \dots , d-1 \} \) such that for all \( i \in \{ 1,
\dots, d-1 \} \), \( x_i \) is a multiple of \( b \) if and only
if \( i \le m \), and the critical points $x'_1, \dots, x'_{d-1}
\in \frak D$ of \( p \) are numbered in such a way that there
exists \( n \in \{0, \dots , d-1 \} \) such that for all \( j \in
\{ 1, \dots, d-1 \} \), \( x'_j \) is a multiple of \( b \) if and
only if \( j \le n \). For every \( k \in \{ 1, \dots, d-1 \} \),
let
\begin{displaymath}
S_k = s_k(x_1, \dots, x_{d-1},0) = s_k(x_1, \dots, x_{d-1}),
\hskip20pt S'_k = s_k(x'_1,\dots,x'_{d-1}).
\end{displaymath}
Then by Equation \eqref{E:3.2},
\begin{equation}
(d-k)S_k = dS'_k, \hskip40pt \forall k \in \{1, \dots, {d-1}\}.
\label{E:6.1}
\end{equation}
When \( k = d-1 \), this gives
\begin{equation}
x_1 \cdots x_{d-1} = b^e x'_1 \cdots x'_{d-1}. \label{E:6.2}
\end{equation}
Since \( b \) is prime in the UFD \( \frak D \), the above
equation implies that \( b \) divides at least one of the roots of
\( p \), that is, \( m \ge 1 \). By definition of \( m \), there
exist \( a_1, \dots, a_m \in \frak D \) such that
\begin{equation}
x_1 = b a_1, \dots, x_m = b a_m. \label{E:6.3}
\end{equation}

\noindent {\rm \textbf{(a)}} Let us prove that \( m = d-1 \).
Assume that \( m < d-1 \). Let \( h \in \mathbb{N}_0 \) be the
highest exponent such that \( b^h \) divides \( (m+1) \) in \(
\frak D \), and let \( r = e - h \). Since \( m+1 < d = b^e \), we
have by Lemma \ref{T:6.2} that \( h < e \), so that \( r>0 \).
Furthermore, by definition of \( h \) there exists \( c \in \frak
D \) such that \( (m+1) = b^h c \) and \( b \) does not divide \(
c \). By applying Equation~\eqref{E:6.1} with $k = d-1-m$, we
obtain \( (m+1) S_k = b^eS'_k \), that is \( b^h c S_k = b^e S'_k
\), that is  \( c S_k = b^r S'_k \), which gives
\begin{align}
c \big[ &x_1 s_{k-1}(x_2, \dots, x_{d-1}) + x_2s_{k-1}(x_3,
\dots, x_{d-1}) + \cdots \notag \\
&+ x_m s_{k-1}(x_{m+1}, \dots, x_{d-1}) + s_k(x_{m+1}, \dots,
x_{d-1}) \big] = b^r S'_k. \label{E:6.4}
\end{align}
Since \( k = d-1-m \), we have \( s_k(x_{m+1}, \dots, x_{d-1}) =
x_{m+1} \cdots x_{d-1} \). Substituting this and \eqref{E:6.3}
into \eqref{E:6.4} gives
\begin{align}
c x_{m+1} \cdots x_{d-1} = b \big[ b^{r-1} S'_k &- c a_1
s_{k-1}(x_2, \dots, x_{d-1}) \notag \\
&{}- \cdots - c a_m s_{k-1}(x_{m+1}, \dots, x_{d-1}) \big]. \notag
\end{align}
Since \( b \) does not divide \( c \), and \( r \) is a positive
integer, it follows that \( b \) divides at least one of the roots
\( x_{m+1}, \dots, x_{d-1} \), which contradicts the definition of
\( m \). Therefore \( m < d-1 \) is impossible, so that we do have
\( m = d-1 \), which implies (a).

\noindent {\rm \textbf{(b)}} If \( n = d-1 \), then (b) is true.
Assume that \( n \le d-2 \). Let \( s = d-1-e \). Let us prove
that \( n \ge s \). Assume that \( n < s = d-1-e \). Then \( d >
n+1+e \ge 0+1+1 =2 \), so that \( d \ge 3 \). Substituting
\eqref{E:6.3} into Equation \eqref{E:6.2}, knowing from (a) that
\( m = d-1 \), gives
\begin{equation}
b^{d-1} a_1 \cdots a_{d-1} = b^e x'_1 \cdots x'_{d-1}, \notag
\end{equation}
that is
\begin{equation}
b^s a_1 \cdots a_{d-1} = x'_1 \cdots x'_{d-1}. \label{E:6.5}
\end{equation}
Let us show that \( s \ge 1 \). If \( e=1 \), then \( s = d-1-e
\ge 3-1-1 = 1 \). If \( e \ge 2 \), then \( s = d-1-e = b^e-1-e
\ge 2^e-1-e \ge 1 \), by induction on \( e \), starting from \( e
= 2 \). Thus \( s \ge 1 \) in all cases, and Equation
\eqref{E:6.5} implies that \( b \) divides at least one of the
critical points of \( p \), that is, \( n \ge 1 \). Then by
definition of \( n \), there exist \( a'_1, \dots, a'_n \in \frak
D \) such that
\begin{equation}
x'_1 = b a'_1, \dots, x'_n = b a'_n. \label{E:6.6}
\end{equation}

Let \( k = d-1-n \) and \( t = k - e \). Since \( n \le d-2 \), we
have that \( k \ge 1 \). Let \( A_k = s_k(a_1, \dots, a_{d-1}) \).
By~\eqref{E:6.1}, we have \( (n+1) S_k = b^e S'_k \), that is \(
(n+1) b^k A_k = b^eS'_k \), that is, \( (n+1) b^t A_k = S'_k \),
which gives
\begin{align}
(n+1) b^t A_k = x'_1 s_{k-1}(x'_2, \dots, x'_{d-1} ) & {}+ \cdots
+ x'_n s_{k-1}( x'_{n+1}, \dots, x'_{d-1} ) \notag\\
& {} + s_k( x'_{n+1}, \dots, x'_{d-1} ). \label{E:6.7}
\end{align}
Since \( k = d-1-n \), we have \( s_k( x'_{n+1}, \dots, x'_{d-1} )
= x'_{n+1} \cdots x'_{d-1} \). Substituting this and \eqref{E:6.6}
into \eqref{E:6.7} gives
\begin{align}
x'_{n+1} \cdots x'_{d-1} = b \big[ (n+1) b^{t-1} A_k &- a'_1
s_{k-1}(x'_2, \dots, x'_{d-1} ) \notag\\
&{}- \cdots - a'_n s_{k-1}( x'_{n+1}, \dots, x'_{d-1} ) \big].
\notag
\end{align}
We have \( t = k-e = d-1-n-e = s - n > 0 \), so that \( t \ge 1
\), and the above equation implies that \( b \) divides at least
one of the critical points \( x'_{n+1}, \dots, x'_{d-1} \). But
this contradicts the definition of \( n \). Therefore, \( n < s \)
is impossible, we do have \( n \ge s = d-1-e \), which implies
(b).

\noindent {\rm \textbf{(c)}} Assume that zero is a root of \( p \)
of multiplicity one, and that \( b \) divides \( e \) in \( \frak
D \). Let us prove (c2) before proving (c1).

\noindent {\rm \textbf{(c2)}} Assume that \( n = d-1-e \). Since
by hypothesis \( b \) divides \( e \) in~\( \frak D \), there
exists \( f \in \frak D \) such that \( e = bf \). Let \( g =
b^{e-1} -f \), and let \( A_e = s_e( a_1, \dots, a_{d-1} ) \).
From \( 2^e \le b^e = d < 2^d  \), we get that \( e \le d-1 \).
Hence, by Equation \eqref{E:6.1}, we have \( (d-e)S_e = d S'_e \),
that is \( (d-e) b^e A_e = b^e S'_e \), that is \( (b^e-bf) A_e =
S'_e \), that is \( bg A_e = S'_e \), that is
\begin{align}
b g A_e = x'_1 s_{e-1}(x'_2, \dots, x'_{d-1} ) & {} +
\cdots + x'_{n} s_{e-1}( x'_{n+1}, \dots, x'_{d-1} ) \notag \\
& {} + s_e( x'_{n+1}, \dots, x'_{d-1} ). \label{E:6.8}
\end{align}
Since  \( e = d-1-n \), we have  \( s_e ( x'_{n+1}, \dots,
x'_{d-1} ) = x'_{n} \cdots x'_{d-1} \). Substituting this and
\eqref{E:6.6} into \eqref{E:6.8} gives
\begin{align}
x'_{n+1} \cdots x'_{d-1} = b \big[g A_e & {} - a'_1
s_{e-1}(x'_2, \dots, x'_{d-1} ) - \cdots \notag \\
& {} - a'_{n} s_{e-1}( x'_{n+1}, \dots, x'_{d-1} ) \big], \notag
\end{align}
which implies that at least one of the critical points \( x'_{n +
1}, \ dots, x'_{d-1} \) is a multiple of \( b \). This contradicts
the definition of \( n \). Therefore \( n = d-1-e \) is
impossible. Since \( n \ge d-1-e \) by (b), it follows that \( n
\ge d-e = (d-1) - (e-1) \), which implies (c2).

\noindent {\rm \textbf{(c1)}} Substituting \eqref{E:6.3} and
\eqref{E:6.6} into Equation \eqref{E:6.2}, knowing from (a) and
(c2) that \( m = d-1 \) and \( n \ge d - e \), gives
\begin{displaymath}
b^{d-1} a_1 \cdots a_{d-1} = b^e b^{d-e} a'_1 \cdots a'_{d-e}
x'_{d-e+1} \cdots x'_{d-1},
\end{displaymath}
that is
\begin{displaymath}
a_1 \cdots a_{d-1} = b a'_1 \cdots a'_{d-e} x'_{d-e+1} \cdots
x'_{d-1}.
\end{displaymath}
This implies that \( b \) divides at least one of the elements \(
a_1 \), \dots, \( a_{d-1} \) of~\( \frak D \), that is, there
exists \( i \in \{1, \dots , d-1 \} \) and \( b_i \in \frak D \)
such that \( a_i = b b_i \). Consequently, \( x_i = b a_i = b^2
b_i \), so that (c1) is true.
\renewcommand{\qedsymbol}{$\blacksquare$}
\end{proof}

\noindent \textbf{Remarks.} \textbf{(1)} In Theorem \eqref{E:6.3},
the hypothesis that at least one root of $p$ is zero is not a loss
of generality, because if $x_0$ denotes a root of a polynomial
$p_1 \in \frak D[x]$, then zero is a root of $p_2(x) = p_1(x +
x_0)$, and by Corollary \ref{T:2.2}, $p_1$ is nice if and only if
$p_2$ is nice.

\noindent \textbf{(2)} If we don't make this shift, then the
conclusion of Theorem \eqref{E:6.3} (a) and (b) is that there
exists at least \( d -1 - e \) critical points of \( p \) such
that if $\mathcal S$ denotes the union of these critical points
with the roots of $p$, then all of the differences of two elements
of $\mathcal S$ are multiples of $d$.

\noindent \textbf{(3)} By Proposition \ref{T:6.1}, a ring of
algebraic integers of an algebraic number field satisfies the
hypotheses of Theorem \ref{T:6.3} iff it is a UFD. For more
information about rings of algebraic integers that are UFDs, see
\cite[pp. 121--123, p. 140, and pp. 185--186]{Narkiewicz}.

\noindent \textbf{(4)} In the case of a computer search for nice
polynomials whose degree~$d$ is a power \( e \) of a prime \( b \)
in $\mathbb{Z}[x]$, we may look for integers \( x_1 \), \dots, \(
x_{d-1} \), \( x_d = 0 \),  \( x'_1 \), \dots, \( x'_{d-1} \)
satisfying \eqref{E:3.2}. Then by Theorem \ref{E:6.3}, we have to
consider only systems such that \( d -1 \) roots and \( d - 1 - e
\) critical points are multiples of $b$, which reduces the number
of choices for each of them by a factor $b$. Thus Theorem
\ref{T:6.3} makes such a search $b^{d-1} \cdot b^{d-1-e} =
b^{2d-2-e}$ times faster.

\begin{corollary}\label{T:6.4}
{\rm [Application to nice quartics].} Let $p \in \mathbb{Z}[x]$ be
a nice quartic having zero as a root of multiplicity one.
Then: \\
\phantom{.} \hskip10pt {\rm \textbf{(a)} }
All of the roots of $p$ are multiples of $2$. \\
\phantom{.} \hskip10pt {\rm \textbf{(b)} } At least one nonzero
root of \( p \) is a multiple of
\( 4 \). \\
\phantom{.} \hskip10pt {\rm \textbf{(c)} } At most \( 1 \) of the
critical points of $p$ is not a multiple of $2$.
\end{corollary}

\begin{proof}
The conclusion is a direct consequence of Theorem \ref{T:6.3}
applied with \( b = e = 2 \).
\renewcommand{\qedsymbol}{$\blacksquare$}
\end{proof}

\textbf{Remark. } Corollary \ref{T:6.4} is not valid if we replace
\( \mathbb Z \) by \( \mathbb G \), because \( 2 = (1 + i)(1 - i)
\) is not prime in \( \mathbb G \).


\section{Nice polynomials with three roots}

\noindent In this section, we present an important reduction of
the system of equations for the roots and critical points of nice
polynomials of arbitrary degree with exactly three roots. This is
done in Theorem 7.1. For the proof of this theorem, it is
essential to extend the definition of the binomial coefficients,
in order to avoid serious problems with finite sums over different
ranges with too many cases to consider. In order to work on a
common range, we have to replace all of theses finite sums by
infinite sums over all of the integers, where all of these
infinite sums have only a finite number of nonzero terms.

\noindent \textbf{Extended binomial coefficients.} Using a small
part of \cite{HP}, we extend the definition of the binomial
coefficients by defining \( \binom{n}{k} = 0 \), for all \( n \in
\mathbb N_0 \) and \( k \in \mathbb Z \) such that \( k < 0 \), or
\( k > n \). By \cite{HP}, these extended binomial coefficients
still satisfy the basic identities:
\begin{equation}
\binom{n-1}{k-1} + \binom{n-1}{k} = \binom{n}{k}, \hskip30pt
\binom{n-1}{k-1} = \frac{k}{n} \binom{n}{k}, \label{E:7.1}
\end{equation}
for all \( n \in \mathbb{N} \), \( k \in \mathbb{Z} \).

In the following theorem, we consider a nice polynomial \( p \)
with three roots \( r_0 \), \( r_1 \), and \( r_2 \) of
multiplicities \( m_0 \), \( m_1 \), and \( m_2 \). By Corollary
\ref{T:2.2}, one of the roots of \( p \), say \( r_0 \), can be
chosen as zero. It is well known that \( p' \) has the same three
roots with multiplicities at least equal to \( m_0-1 \), \( m_1-1
\), and \( m_2-1 \). Since the degree of \( p' \) is \( d-1 \), \(
p' \) must have two additional roots \( r'_1 \) and \( r'_2 \),
not necessarily distinct from \( r_0 = 0 \), \( r_1 \), and \( r_2
\). By Corollary \ref{T:3.2}, the \( (m_1 + m_2) \) nonzero roots
\( r_1 \) and \( r_2 \) of \( p \), and the \( (m_1 + m_2) \)
critical points \( r_1 \), \( r_2 \), \( r'_1 \), and \( r'_2 \)
of \( p \), counted with their multiplicities, satisfy a system of
\( (m_1 + m_2) \) equations of degrees \( 1 \), \( 2 \), \dots, \(
m_1 + m_2 \). Theorem \ref{T:7.1} shows that we can reduce this
system to a system of two equations of degrees one and two for the
four variables \( r_1 \), \( r_2 \), \( r'_1 \), and \( r'_2 \).

\begin{theorem}\label{T:7.1}
{\rm [Reduction of the system of equations \eqref{E:3.2} for nice
polynomials with three roots to a system of two equations of
degrees one and two.]} Let \( \frak D \) be an integral domain of
characteristic zero. Let \( r_1 \) and \( r_2 \) be nonzero
elements of \( \frak D \). Let \( r'_1 \) and \( r'_2 \) be
elements of \( \frak D \). Let \( m_0, m_1, m_2 \in \mathbb N \).
Let \( d = m_0 + m_1 + m_2 \). Let
\begin{align}
p(x) &= x^{m_0}(x-r_1)^{m_1}(x-r_2)^{m_2}, \notag \\
q(x) &= d x^{m_0-1}(x-r_1)^{m_1-1}(x-r_2)^{m_2-1}(x-r'_1)(x-r'_2).
\notag
\end{align}
Then \( q = p' \) iff
\begin{align}
(d-m_1)r_1 + (d-m_2)r_2 &= d(r'_1 + r'_2), \label{E:7.2} \\
m_0 r_1 r_2 &= d r'_1 r'_2. \label{E:7.3}
\end{align}
\end{theorem}

\begin{proof}
For every \( i \in \{1, \dots, d \} \), let
\begin{displaymath}
x_i =
\begin{cases}
r_1, &\text{if } i \le m_1, \\
r_2, &\text{if } m_1+1 \le i \le m_1 + m_2, \\
0, &\text{if } i \ge m_1 + m_2 + 1.
\end{cases}
\end{displaymath}
For every \( j \in \{1, \dots, d-1 \} \), let
\begin{displaymath}
x'_j =
\begin{cases}
r_1, &\text{if } j \le m_1-1, \\
r_2, &\text{if } m_1 \le j \le m_1 + m_2 - 2, \\
r'_1, &\text{if } j = m_1 + m_2 - 1, \\
r'_2, &\text{if } j = m_1 + m_2, \\
0, &\text{if } j \ge m_1 + m_2 + 1.
\end{cases}
\end{displaymath}

\noindent For every $k \in \{1, \dots, d-1\}$, let
\begin{alignat}{2}
S_k &= s_k(x_1, \dots, x_d) &&= s_k(x_1, \dots, x_{m_1 + m_2} ), \notag \\
S'_k &= s_k(x'_1, \dots, x'_{d-1}) &&= s_k(x'_1, \dots, x'_{m_1 +
m_2}), \notag
\end{alignat}
\begin{displaymath}
\mathcal{E}(k) = \big[ (d-k) S_k = d S'_k \big]. \notag
\end{displaymath}
Since \( S_k = S'_k = 0 \) for all \( k \in \{1, \dots, d-1\} \)
such that \( k > m_1 + m_2 \), we have by Corollary \ref{T:3.2}
that
\begin{equation}
(q = p') \iff \big[ \mathcal{E}(k), \hskip10pt \forall k \in \{ 1,
\dots, m_1 + m_2 \} \big]. \label{E:7.4}
\end{equation}
Since
\begin{alignat}{2}
S_1 &= x_1 + \cdots  + x_{m_1+m_2} &&= m_1 r_1 + m_2 r_2, \notag\\
S'_1 &= x'_1 + \cdots  + x'_{m_1+m_2} &&= (m_1-1) r_1 + (m_2-1)
r_2 + r'_1 + r'_2, \notag
\end{alignat}
we have that \( \mathcal{E}(1) \) is equivalent to
\begin{equation}
(d-1)(m_1 r_1 + m_2 r_2) = d \big[(m_1-1) r_1 + (m_2-1) r_2 + r'_1
+ r'_2 \big], \notag
\end{equation}
which is equivalent to \eqref{E:7.2}. Since \( d-m_1-m_2 = m_0 \),
and
\begin{alignat}{2}
S_{m_1+m_2} &= x_1 \cdots x_{m_1+m_2} &&= r_1^{m_1} r_2^{m_2}, \notag \\
S'_{m_1+m_2} &= x'_1 \cdots x'_{m_1+m_2} &&= r_1^{m_1-1}
r_2^{m_2-1} r'_1 r'_2, \notag
\end{alignat}
we have that \( \mathcal{E}(m_1+m_2) \) is equivalent to
\begin{equation}
m_0 r_1^{m_1} r_2^{m_2} = d r_1^{m_1-1} r_2^{m_2-1} r'_1 r'_2,
\notag
\end{equation}
which is equivalent to \eqref{E:7.3}. Thus
\begin{displaymath}
\mathcal{E}(1) \iff \eqref{E:7.2} \hskip15pt \mbox{and} \hskip15pt
\mathcal{E}(m_1 + m_2) \iff \eqref{E:7.3}.
\end{displaymath}
This implies by \eqref{E:7.4} that the proof is complete if  \(
m_1 + m_2 = 2 \), and otherwise it remains to prove that
\eqref{E:7.2} and \eqref{E:7.3} imply \( \mathcal{E}(k) \) for all
\( k \in \{ 2, \dots, m_1 + m_2 -1 \} \).

Assume that \( m_1 + m_2 > 2 \), and that \eqref{E:7.2} and
\eqref{E:7.3} hold. Let \( k \in \{ 2, \dots, m_1 + m_2 -1 \} \).
Let us prove that \( \mathcal{E}(k) \) is true. With the above
extended definition of the binomial coefficients, we have
\begin{align}
dS'_k = {} &d \sum_{i \in \mathbb Z} \binom{m_1-1}{i}
\binom{m_2-1}{k-i} r_1^i r_2^{k-i} \notag \\
&+ d \sum_{i \in \mathbb Z} \binom{m_1-1}{i}
\binom{m_2-1}{k-1-i} r_1^i r_2^{k-1-i} (r'_1 + r'_2) \notag \\
&+ d \sum_{i \in \mathbb Z} \binom{m_1-1}{i} \binom{m_2-1}{k-2-i}
r_1^i r_2^{k-2-i} r'_1 r'_2, \notag
\end{align}
which gives by \eqref{E:7.2} and \eqref{E:7.3}
\begin{multline}
dS'_k = d \sum_{i \in \mathbb Z} \binom{m_1-1}{i}
\binom{m_2-1}{k-i} r_1^i r_2^{k-i} \notag \\
+ \sum_{i \in \mathbb Z} \binom{m_1-1}{i} \binom{m_2-1}{k-1-i}
r_1^i r_2^{k-1-i} \big[(d-m_1) r_1
+ (d-m_2) r_2\big] \notag \\
+ m_0 \sum_{i \in \mathbb Z} \binom{m_1-1}{i} \binom{m_2-1}{k-2-i}
r_1^i r_2^{k-2-i} r_1 r_2. \notag
\end{multline}
Using \( m_0 = d - m_1 - m_2 \), grouping terms with factor \( d
\), \( m_1 \), and \( m_2 \), simplifying with \eqref{E:7.1} gives
\begin{alignat}{2}
dS'_k &= {}&& d \sum_{i \in \mathbb Z} \binom{m_1}{i}
\binom{m_2}{k-i} r_1^i r_2^{k-i} \notag \\
& && - m_1 \sum_{i \in \mathbb Z} \frac{i}{m_1} \binom{m_1}{i}
\binom{m_2}{k-i} r_1^i r_2^{k-i}  \notag \\
& && -m_2 \sum_{i \in \mathbb Z} \binom{m_1}{i}
\frac{k-i}{m_2} \binom{m_2}{k-i} r_1^i r_2^{k-i}. \notag \\
& {} = {} && (d-k) \sum_{i \in \mathbb Z} \binom{m_1}{i}
\binom{m_2}{k-i} r_1^i r_2^{k-i} = (d-k) S_k. \notag
\end{alignat}
Thus  \( \mathcal{E}(k) \) is true, and the theorem is proved.
\renewcommand{\qedsymbol}{$\blacksquare$}
\end{proof}

\noindent \textbf{Applications. 1.} While the three root case is
not solved yet, Theorem~\ref{T:7.1} makes it easy to find examples
of nice polynomials with three roots of any given multiplicities
with the help of a computer.

\noindent \textbf{2.} We have established the formulas A2, A5, A6,
A7, A8, and A10 of Section 8 below by applying Theorem
\ref{T:7.1}.

\noindent \textbf{Conjecture.} Theorem \ref{T:7.1} suggests that
nice polynomials with one root at zero and \( n \) nonzero roots
may be determined by a reduced system of \( n \) equations of
degrees \( 1 \), \( 2 \), \dots, \( n \) for the \( n \) nonzero
roots and the \( n \) nonzero critical points.


\section{Examples and open problems}

\enlargethispage{10pt}

\noindent In this section, we present examples of nice polynomials
whose coefficients, roots, and critical points are rational
integer or Gaussian integers, with special interest in the case
where the coefficients are real, while the roots and critical
points can be real or complex. For each type of nice polynomials,
we have tried to find the \emph{smallest} one, and to present it
in \emph{its most reduced form}, except for the symmetric ones.
All of the examples of this section have been checked with Maple
through Maple methods totally different from the methods of this
paper.


\noindent \textbf{Definitions.} As usual, we call \emph{diameter}
of a compact subset~\( \mathcal S \) of the complex plane the
largest distance between two points of~\( \mathcal S \). We define
the diameter \( \operatorname{diam}(p) \) of a polynomial \( p \in
\mathbb{G}[x] \) as the diameter of the set of the roots of \( p
\) in the complex plane. We say that a polynomial \( p \in
\mathbb{G}[x] \) is \emph{smaller than} a polynomial  \( q \in
\mathbb{G}[x] \) iff \( \operatorname{diam}(p) <
\operatorname{diam}(q) \). We say that a nice polynomial \( p \in
\mathbb G \) is in \emph{reduced form} iff \( p \) cannot be
transformed into a smaller nice polynomial by using the
transformations of Proposition \ref{T:2.1}. We say that a nice
polynomial \( p \in \mathbb{G} \) is in \emph{most reduced form}
iff  \( p \) is in reduced form, zero is a root of \( p \), \( p
\) cannot be transformed into a nice polynomial \( q \) with  \(
\operatorname{diam}(p) = \operatorname{diam}(q) \) and \( |
\mathcal A(q) | < | \mathcal A(p) | \) by using the
transformations of Proposition \ref{T:2.1}, and when all of the
roots of \( p \) are real, none of them is negative. For example,
\( p(x) = x^2(x-7)^5 \) is a nice polynomial in reduced form, and
\( q(x) = x^5(x-7)^2 \) is the most reduced form of \( p \).

\noindent \textbf{A. A few general formulas.} Here, we present
general formulas that are not in final form. These formulas may
give several times the \emph{same} polynomials, and they do not
separate the cases where some roots are equal. In the following
formulas, \( m \) and \( n \) are rational integers with no common
factors.

\noindent \textbf{A 1.} For every integral domain \( \frak{D} \)
of characteristic zero, and for every positive integer \( n \)
greater than one, the polynomial \( p(x) = x^n \) is a totally
nice polynomial over \( \frak{D} \). It is symmetric iff \( n \)
is even and antisymmetric
iff \( n \) is odd.\\
\textbf{A 2.} For every positive integer \( n \), there exists a
nice antisymmetric polynomial \( p \in \mathbb{Z}[x] \) of degree
\( (2n+1)^2 \) with three roots, given by
\begin{align}
p( x) &= x \big[ x^2 - (2n+1)^2 \big]^{2n(n+1)}, \notag \\
p'(x) &= (2n+1)^2 (x^2-1) \big[ x^2 - (2n+1)^2 \big]^{2n(n+1)-1}.
\notag
\end{align}
\textbf{A 3.} For every positive integer  \( n \), there exists a
nice symmetric polynomial \( p \in \mathbb{Z}[x] \) of degree \(
(2n+1)^2 + 1 \) with two roots, whose derivative is also nice,
which is given by
\begin{align} \allowdisplaybreaks
p( x) &= \big[ x^2 - (2n+1)^2 \big]^{2n(n+1)+1}, \notag \\
p'(x) &= 2 (2n^2 + 2n +1) x \big[ x^2 - (2n+1)^2 \big]^{2n(n+1)},
\notag \\
p''(x) &= 2 (2n+1)^2 (2n^2 + 2n + 1) (x^2-1) \big[ x^2 - (2n+1)^2
\big]^{2n(n+1)-1}. \notag
\end{align}
\textbf{A 4.} Nice polynomials with two roots in reduced form are
given by
\begin{align}
p( x) &= x^{ab}( x - c)^{b(c-a)}, \notag \\
p'(x) &= bc x^{ab-1}( x - c)^{b(c-a)-1}(x-a), \notag
\end{align}
where \( a, b, c \in \mathbb N \) are such that \( c > a \), and
\( \gcd(a,c) = 1 \). \\
\textbf{A 5.} Nice Cubics with three real roots: This case was
solved by M. Chapple, \cite{Chapple} in 1960, Karl Zuser,
\cite{Zuser} in 1963, and Tom Bruggeman and Tom Gush \cite{BG} in
1980. Those cubics are given by the formula
\begin{align}
p( x) &= x \big[ x - 3 m(m + 2 n) \big] \big[ x - 3 n(2 m + n) \big],
\notag \\
p'(x) &= 3( x - 3 mn) \big[ x - (m + 2 n)(2 m + n) \big]. \notag
\end{align}
Karl Zuser observed that when \( mn \equiv 1 \pmod 3 \), these
polynomials can be made three times smaller by using
Proposition~\ref{T:2.1}(c). \\
\textbf{A 6.} Nice Cubics with one real root, and two complex
conjugate roots with even imaginary part are given by
\begin{align}
p( x) &= ( x^2 + 36 m^2n^2) \big[ x - 3( m^2 + 3 n^2) \big], \notag \\
p'( x) &= 3( x - 2 m^2)( x - 6 n^2). \notag
\end{align}
\textbf{A 7.} Nice Cubics with one real root, and two complex
conjugate roots with odd imaginary part are given by
\begin{align}
p( x) &= \big[ x^2 + 9 (2m + 1)^2(2n + 1)^2 \big] \notag \\
& \phantom{8888888888888888} \cdot \Big[ x - 6
\big[ m(m + 1) + 3n( n + 1) + 1 \big] \Big], \notag \\
p'( x) &= 3[ x - (2m+1)^2][ x - 3 (2n + 1)^2]. \notag
\end{align}
\textbf{A 8.} Nice quartics with a double root are given by
\begin{align}
p( x) &= x^2 \big[ x - 2(2m^2 -mn - n^2) \big] \big[ x - 2(2m^2 +
mn - n^2 ) \big], \notag \\
p'( x) &= 4 x \big[ x - 2(m^2 - n^2) \big] \big[ x - (4 m^2 - n^2)
\big], \notag
\end{align}
and by
\begin{align}
p( x) &= x \big[ x - 2(n^2 + mn - 2m^2) \big]^2 ( x - 4mn ),
\notag \\
p'( x) &= 4 \big[ x - (2mn + n^2) \big] \big[ x - 2(n^2 + mn -
2m^2) \big] \big[ x - 2m(n - m) \big]. \notag
\end{align}
\textbf{A 9.} Nice symmetric quartics with four real roots. This
case was solved by Chris Caldwell \cite{Caldwell} in 1990. These
quartics are given by
\begin{align}
p( x) &= \big[ x^2 - (n^2-m^2 - 2 m n )^2 \big] \big[ x^2
- (n^2-m^2 + 2 m n )^2 \big], \notag \\
p'( x) &= 4 x \big[ x^2 - (m^2 + n^2)^2 \big]. \notag
\end{align}
\textbf{A 10.} Nice quintics with a triple root at zero, and a
critical point that is a multiple of 15 are given by
\begin{align}
p( x) &= x^3 \big[ x - 5m (m - 4 n ) \big]
\big[ x - 5n (4m - 15 n ) \big], \notag \\
p'( x) &= 5 x^2 \big[ x - 15n(m - 4 n) \big] \big[ x - m(4 m -15
n) \big]. \notag
\end{align}
Other nice quintics with a triple root at zero are given by
\begin{align}
p( x) &= x^3 \big[ x - 5m (4 n - 5 m ) \big]
\big[ x - 5 n (3n - 4 m ) \big], \notag \\
p'( x) &= 5 x^2 \big[ x - 5m (3n - 4m) \big] \big[ x - 3n(4n - 5m)
\big]. \notag
\end{align}
\textbf{A 11.} Nice symmetric sextics with at least two double
roots are given by
\begin{align}
p( x) &= \big[ x^2 - (m^2 + 4 m n - 2 n^2)^2 \big]^2
\big[ x^2 - (m^2 - 2 m n - 2 n^2)^2 \big], \notag \\
p'( x) &= 6 x \big[ x^2 - (2 n^2 - 4 m n -m^2)^2 \big] \big[ x^2 -
(m^2 + 2 n^2)^2 \big]. \notag
\end{align}

\noindent \textbf{B. A few numerical examples.}

\noindent \textbf{B 1.} The smallest nice cubic with three
distinct real roots:
\begin{displaymath}
p( x) = x (x - 9)(x - 24), \hskip30pt p'( x) = 3(x-4) (x-18).
\end{displaymath}
\noindent \textbf{B 2.} The smallest nice cubic with two complex
conjugate roots:
\begin{displaymath}
p( x) = (x - 6)(x^2 + 9), \hskip30pt p'( x) = 3(x-1) (x-3).
\end{displaymath}
\noindent \textbf{B 3.} The smallest nice cubic with roots at
unequal distances from each other in the complex plane:
\begin{displaymath}
p( x) = x \big[ x - 3(1 + 4i) \big](x-15), \hskip5pt p'( x) =
3\big[ x - 3(1 + 2i) \big] \big[ x - (9 + 2i) \big].
\end{displaymath}
\textbf{B 4.} The smallest nice symmetric quartic with four
distinct complex roots is \(p( x) = x^4 - 1 \), which gives \( p'(
x) = 4 x^3 \).\\
\textbf{B 5.} The smallest nice symmetric quartic with four
distinct real roots:
\begin{displaymath}
p( x) = (x^2-1) (x^2-49), \hskip30pt p'( x) = 4 x (x^2-25).
\end{displaymath}
\noindent \textbf{B 6.} Concerning the case of nice non-symmetric
quartics with four distinct real roots, the first five examples
were published by Chris Caldwell \cite{Caldwell} in 1990. Since
then, it was not known whether other examples existed. After
considerable simplification of the computation with the results of
this paper, and formulas that are not in publishable form yet, we
have found 358 other examples with the help of a computer. The
five examples of \cite{Caldwell} are in the positions 88, 107,
181, 247, and 321 on our list of examples ordered in increasing
order. The smallest of our examples is:
\begin{align}
p( x) &= x(x-50)(x-176)(x-330), \notag \\
p'( x) &= 4 (x-22)(x-120)(x-275). \notag
\end{align}
The largest example that we have obtained, in position 363 on our
list, after two weeks of uninterrupted computation with a
computer, is
\begin{align}
p( x) &= x (x - 33408) (x - 44100) (x - 138040), \notag \\
p'( x) &= 4 (x - 11760) (x - 38976) (x - 110925). \notag
\end{align}
\textbf{B 7.} The only nice quintic with five distinct real roots
that we have found is
\begin{align}
p( x) &= x(x-180)(x-285)(x-460)(x-780), \notag \\
p'( x) &= 5(x-60)(x-230)(x-390)(x-684). \notag
\end{align}
\textbf{B 8.} We have found only two nice quintics with five
distinct complex roots. The first one has a double critical point:
\begin{align}
p( x) &= x (x - 595) (x - 1020) \big[ (x -220)^2 + 40,000 \big], \notag \\
p'( x) &= 5  (x - 170)^2(x - 420) (x - 884), \notag \\
p( x) &= x (x - 585) (x - 1040) \big[ (x - 270)^2 + 44,100 \big], \notag \\
p'( x) &= 5 (x - 130)(x - 312) (x - 390)(x-900). \notag
\end{align}
\textbf{B 9.} The smallest nice antisymmetric polynomial with more
than one root is \( p(x) = x(x^2-9)^4 \), with \( p'(x) = 9
(x^2-1) (x^2-9)^3 \). \\
\textbf{B 10.} The smallest nice symmetric polynomial with more
than one critical point, and whose first derivative is also nice,
is:
\begin{displaymath}
p(x) = (x^2-9)^5, \hskip5pt p'(x) = 10 x (x^2-9)^4, \hskip5pt
p''(x) = 90 (x^2-1) (x^2-9)^3.
\end{displaymath}

\vskip10pt

\noindent \textbf{Open problems.}

\noindent \textbf{1.} With the equivalence relation defined in
Section 2, find the number of equivalence classes for any set of
nice polynomials of a given type (For example, symmetric with four
roots of any multiplicities). \\
\textbf{2.} Generalize Theorem \ref{T:6.3} from polynomials whose
degree is a prime power to polynomials of any degree. A first step
would be to consider polynomials whose degree is the product of
two primes. \\
\textbf{3.} Prove or disprove the conjecture of the end of Section
7. \\
\textbf{4.} Find a nice sextic with six distinct roots. \\
\textbf{5.} Find a nice antisymmetric polynomial with more than
three roots. \\
\textbf{6.} Find a general formula for nice polynomials with three
roots by using the system of two equations \eqref{E:7.2} and
\eqref{E:7.3}. \\
\textbf{7.} Find a general formula for nice non-symmetric
quartics.

\noindent \textbf{Acknowledgment.} The author would like to
express his gratitude to the University of Toledo for its strong
support and encouragements of this project during his visiting
position there from 2000 to 2003.


\enlargethispage{5pt}

\vskip20pt

\noindent Jean-Claude Evard \\
Department of Mathematics \\
Western Kentucky University \\
Bowling Green, KY 42101--3576

\noindent E-mail: Jean-Claude.Evard@wku.edu

\end{document}